\documentclass[letterpaper, 10 pt, conference]{ieeeconf}  
\IEEEoverridecommandlockouts                              
\usepackage[Symbol]{upgreek}
\usepackage{amsmath,amssymb,amsfonts,mathrsfs,mathtools,dsfont,mathrsfs,bm,bbm} 
\DeclareMathAlphabet{\mathpzc}{OT1}{pzc}{m}{it}
\DeclareSymbolFontAlphabet{\amsmathbb}{AMSb}%
\usepackage{csvsimple,booktabs}
\usepackage{float}
\usepackage[caption=false]{subfig}
\usepackage{graphicx,epsfig,xcolor}
\usepackage{balance}

\let\labelindent\relax
\usepackage{enumitem}

\usepackage{cite}
\makeatletter
\newcommand{\removelatexerror}{\let\@latex@error\@gobble}
\let\NAT@parse\undefined
\makeatother
\usepackage[hyphens]{url}
\usepackage[breaklinks, colorlinks, linkcolor=purple, anchorcolor=purple, citecolor=purple, urlcolor=purple]{hyperref}

\usepackage[ruled,vlined,linesnumbered]{algorithm2e}

\usepackage{fancyhdr} 

\newcommand{\beq}{\begin{equation}}
\newcommand{\eeq}{\end{equation}}
\newcommand{\beqa}{\begin{eqnarray}}
\newcommand{\eeqa}{\end{eqnarray}}
\newcommand{\beqan}{\begin{eqnarray*}}
\newcommand{\eeqan}{\end{eqnarray*}}

\newcommand\T{{\mathpalette\raiseT\intercal}}
\newcommand\raiseT[2]{\raisebox{0.25ex}{$#1#2$}
}

\renewcommand{\dim}{{\sf dim}} 

\newcommand{\Aset}{\mathds{A}}

\newcommand{\Eset}{\mathds{E}}

\newcommand{\Rset}{\mathds{R}}

\newcommand{\Xset}{\mathds{X}}

\newcommand{\Acal}{{\cal A}}

\newcommand{\Ical}{{\cal I}}
\newcommand{\Jcal}{{\cal J}}

\newcommand{\Ncal}{{\cal N}}

\newcommand{\argmin}{\mathop{\rm argmin}}

\newcommand{\bone}{\mathbf{1}}

\renewcommand{\v}[1]{{\bm{#1}}}

\newcounter{l1}
\newcounter{l2}
\newcounter{l3}
\newcommand{\bdotlist}{\begin{list}{$\bullet$}{}}
\newcommand{\bboxlist}{\begin{list}{$\Box$}{}}
\newcommand{\bbboxlist}{\begin{list}{\raisebox{.005in}{{\tiny
$\blacksquare$ \ \ }}}{}}
\newcommand{\bdashlist}{\begin{list}{$-$}{} }
\newcommand{\blist}{\begin{list}{}{} }
\newcommand{\barablist}{\begin{list}{\arabic{l1}}{\usecounter{l1}}}
\newcommand{\balphlist}{\begin{list}{(\alph{l2})}{\usecounter{l2}}}
\newcommand{\bAlphlist}{\begin{list}{\Alph{l2}.}{\usecounter{l2}}}
\newcommand{\bdiamlist}{\begin{list}{$\diamond$}{}}
\newcommand{\bromalist}{\begin{list}{(\roman{l3})}{\usecounter{l3}}}

\newtheorem{theorem}{Theorem}
\newtheorem{lemma}{Lemma}
\newtheorem{proposition}{Proposition}
\newtheorem{corollary}{Corollary}
\newtheorem{remark}{Remark}

\newtheorem{assumption}{Assumption}

\renewcommand{\bone}{{\mathds{1}}}

\newcommand{\nbd}{{\sf nbd}}

\definecolor{teal}{rgb}{0.0, 0.5, 0.5}

\begin{document}

\title{\LARGE {\bf{Distributed Multi-Area Optimal Power Flow \\ via Rotated Coordinate Descent Critical Region Exploration}}}

\author{Haitian Liu \qquad Ye Guo \qquad Hongbin Sun  \qquad Weisi Deng
\thanks{H. Liu and Y. Guo are with the Tsinghua-Berkeley Shenzhen Institute, Shenzhen, Guangdong 518055, China. H. Sun is with the Tsinghua University, Beijing 100084, China. W. Deng is with the China Southern Power Grid, Guangzhou, Guangdong 510530, China. 
This work is supported in part by the National Science Foundation of China under Grant 51977115.}
\thanks{Corresponding author: Ye Guo, e-mail: guo-ye@sz.tsinghua.edu.cn. }
} 

\maketitle
\thispagestyle{fancy} 
\lhead{} 
\chead{} 
\rhead{} 
\lfoot{} 
\cfoot{} 
\rhead{\thepage} 
\renewcommand{\headrulewidth}{0pt} 
\renewcommand{\footrulewidth}{0pt} 
\pagestyle{fancy}  
\rhead{\thepage}  

\begin{abstract}
  We consider the problem of distributed optimal power flow (OPF) for multi-area electric power systems. A novel distributed algorithm is proposed, referred to as the rotated coordinate descent critical region exploration (RCDCRE). It allows each entity to independently update its boundary information and optimally solve its local OPF in an asynchronous fashion. RCDCRE method stitches coordinate descent and parametric programming using coordinate system rotation to reduce coordination, keep privacy and ensure convergence. The solution process does not require warm starts and can iterate from infeasible initial points using penalty-based formulations. The effectiveness of RCDCRE is verified based on IEEE 2-area 44-bus and 4-area 472-bus systems.
\end{abstract}




\section{Introduction}
\label{sec:intro}

Decarbonization is a common consensus to hedge against global warming, which pushes the need to increase the installed capacity and utilization level of renewable energy resources (RESs). 
As the spatial distribution of RESs is unbalanced, some areas cannot satisfy local demand unless a substantial amount of power is imported from neighboring areas. In contrast, some other areas lack enough flexibility to accommodate the RESs, leading to significant curtailment.
Due to privacy and computation capability concerns, the modern vast interconnected power system is jointly operated by multiple utilities. 
Each entity is responsible for optimizing the dispatch of local RESs according to the interface power exchange. 
Hence, optimal tie-line scheduling is crucial for decarbonization. 
Moreover, it also contributes to the system's overall efficiency and provides resilience to an area under extreme weather conditions.

As gathering all data at a central location may not be the best, a distributed computational paradigm is advocated to coordinate multi-area systems. 
The utilities collectively solve a distributed OPF to obtain the interchange schemes by optimizing internal assets and exchanging intermediate variables. 
Current distributed methods tailored for the joint power system operation are mainly derived from dual decomposition. See \cite{molzahnSurveyDistributedOptimization2017, kargarianDistributedDecentralizedDC2018, patariDistributedOptimizationDistribution2021} for a comprehensive review.
Dual decomposition enjoys a near block-separable structure by adopting the (augmented) Lagrangian relaxation (LR) techniques.
The distributed multi-area OPF can be solved by iteratively updating the local dispatch orders and coupling constraints' dual multipliers.
The dual update also makes economic sense, which indicates the price signals of net power exchange and the shadow prices of tie-line capacity. 

However, since coupling constraints are relaxed, all intermediate solutions are not feasible until asymptotic convergence.
Dual methods are generally first-order methods that optimize over a non-smooth dual function. They may suffer from slow convergence in large multi-area power systems.
In practice, an infeasible dispatch may be returned as the final solution due to the limited computation time. 
Infeasible dispatch is less desirable for two reasons, 1) The economic performance is degraded if multiple automatic generation control units are activated to make up for the power mismatch. 2) Or even worse, the returned solutions cannot dispatch, causing security issues.
Besides, the LR decomposition has strict assumptions on the problem structures to ensure feasible primal solutions can be obtained in the limit\cite{necoaraLinearConvergenceDistributed2015a}.
Divergent cases may occur when applying the alternating direction method of multipliers (ADMM) to a system more than 2-area\cite{chenDirectExtensionADMM2016}.

Meanwhile, primal decomposition has also made progress in solving the multi-area OPF problem. Typical methods includes marginal equivalent decomposition (MED)\cite{zhaoMarginalEquivalentDecomposition2014}, column and constraint generation (C\&CG)\cite{zengSolvingTwostageRobust2013}, and critical region exploration (CRE)\cite{guoRobustTieLineScheduling2018, guoCoordinatedMultiAreaEconomic2017}. 
MED identifies the marginal variables to exploit the finite structure of active constraint sets. 
C\&CG improves the convergence of Benders decomposition by simultaneously submitting cutting planes and binding variables.
CRE leverages multi-parameter programming (mpP) theory, finding optima with the finite search of critical regions. 
Empirically, primal decomposition methods have better feasibility guarantees and convergence over the duals, as there is no constraint relaxed. 
However, some methods may exchange actual costs and network information (e.g., MED, C\&CG), suffer from the strong coordination with synchronization requirements (e.g., C\&CG, CRE), or unable to dispatch when current boundary states are infeasible to the local problems (e.g., MED, CRE).

To date, solving the large-scale distributed multi-area OPF problem is still a challenging problem, both from the primal and dual decomposition. 
This paper proposed a \emph{rotated coordinate descent critical region exploration} (RCDCRE) algorithm in terms of primal decomposition and multi-parameter programming theory.
Two merits make RCDCRE tailored to multi-area OPF. 
1) Compared to existing primal decompositions, RCDCRE preserves privacy by exchanging parametric mappings and reduces coordination by separating boundary jurisdiction with an asynchronous update. 
2) Compared to existing dual decompositions, the proposed method can obtain feasible intermediate solutions quickly by parameterizing the penalized intra-regional OPF. Moreover, it also enjoys a faster convergence as more useful information is exchanged. 
Hence, RCDCRE is suitable for the solution process with an infeasible cold start or online implementation with fixed iteration requirements.
Numerical results are encouraging and show convergence superiority over ADMM and Benders method. 



\section{Problem Setup}
\label{sec:model}

\subsection{Multi-Area OPF formulation}

\begin{figure}[htbp]
    \centering
    \vspace{-10pt}
    \includegraphics[width=0.40\textwidth]{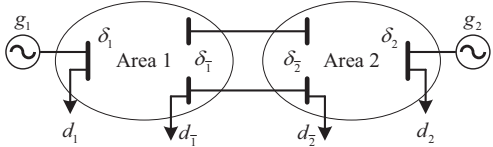}
	\vspace{-10pt}
    \caption{An illustration for multi-area power system.}
    \label{fig:maopf}
\end{figure}

Without loss of generality, we use a two-area network to illustrate the interconnected power system in Figure \ref{fig:maopf}. 

With the DC power flow model for the regional transmission network, an $N$-area OPF problem is formulated as
\begin{subequations}
	\label{eq:maopf_theta}
	\begin{alignat}{2}
	\underset{\v{g}_i, \v{\delta}_i, \v{\delta}_{\bar{i}}}{\text{minimize}} & 
	\quad \sum_{i=1}^{N} \left ( \frac{1}{2} \v{g}_i^\T \v{Q}_i \v{g}_i + \v{c}_i^\T \v{g}_i \right ), 
	\label{eq:maopf_theta.obj} \\ 
	\text{subject to} &
	\quad \v{B}_{i,i} \v{\delta}_i + \v{B}_{i,\bar{i}} \v{\delta}_{\bar{i}} = \v{g}_i - \v{d}_i, 
	\label{eq:maopf_theta.power.int} \\
	&
	\quad \v{B}_{\bar{i},i} \v{\delta}_i + \v{B}_{\bar{i},\bar{i}} \v{\delta}_{\bar{i}} + \sum_{j \in \nbd(i)} \v{B}_{\bar{i},\bar{j}} \v{\delta}_{\bar{j}}
	= - \v{d}_{\bar{i}}, 
	\label{eq:maopf_theta.power.bnd} \\
	&
	\quad - \v{f}_i \leq \v{H}_{i,i} \v{\delta}_i + \v{H}_{i,\bar{i}} \v{\delta}_{\bar{i}} \leq \v{f}_i,
	\label{eq:maopf_theta.line.int} \\
	&
	\quad - \v{f}_{\bar{i}} \leq \v{H}_{\bar{i},\bar{i}} \v{\delta}_{\bar{i}} + \sum_{j \in \nbd(i)} \v{H}_{\bar{i},\bar{j}} \v{\delta}_{\bar{j}} \leq \v{f}_{\bar{i}}, 
	\label{eq:maopf_theta.line.TL} \\	
	&
	\quad \underline{\v{g}}_i \leq \v{g}_i \leq \overline{\v{g}}_i, 
	\label{eq:maopf_theta.gen.int} \\
	&
	\quad - \v{\pi} \leq \v{\delta}_i  \leq \v{\pi}, 
	\label{eq:maopf_theta.angle.int} \\
	&
	\quad - \v{\pi} \leq \v{\delta}_{\bar{i}} \leq \v{\pi}, \quad \forall i = 1 , \dots, N 
	\label{eq:maopf_theta.angle.bnd} \\	
    & \quad \delta^{\text{ref}} = 0. \label{eq:maopf_theta.ref}
	\end{alignat}
\end{subequations}
where, as shown in Figure \ref{fig:maopf}, decision variables include area $i$'s power generation $\v{g}_i$, internal and boundary phase angles $\v{\delta}_i$, $\v{\delta}_{\bar{i}}$. 
The objective \eqref{eq:maopf_theta.obj} is to minimize the sum of each area's generation costs.
The DC power flow equation is divided into constraints \eqref{eq:maopf_theta.power.int}-\eqref{eq:maopf_theta.power.bnd}. 
Notation $\nbd(i)$ collects all adjacent areas of area $i$. 
Constraints \eqref{eq:maopf_theta.line.int}-\eqref{eq:maopf_theta.line.TL} describe internal and tie-line flow limits respectively. 
The lower and upper bounds of generation capacities, internal and boundary phase angle limits are summarized in constraints \eqref{eq:maopf_theta.gen.int}-\eqref{eq:maopf_theta.angle.bnd}. 
We assume there is no generator on the boundary of each area. 
Constraint \eqref{eq:maopf_theta.ref} artificially assigns a reference phase angle. 
To simplify notations, we define new decision variables $\v{x}_i$ and $\v{\theta}_i$ as 
\begin{alignat*}{2}
	\v{x}_i = \left [ \v{g}_i^\T,\ \v{\delta}_i^\T \right ]^\T, \quad \v{\theta}_i = \v{\delta}_{\bar{i}}
\end{alignat*}
By summarizing \eqref{eq:maopf_theta.power.int}, \eqref{eq:maopf_theta.power.bnd}, \eqref{eq:maopf_theta.line.int}, \eqref{eq:maopf_theta.gen.int}, and \eqref{eq:maopf_theta.angle.int} into \eqref{eq:probP.local}, while \eqref{eq:maopf_theta.line.TL} and \eqref{eq:maopf_theta.angle.bnd} become instances of \eqref{eq:probP.couple}. 
We can also eliminate internal phase angles $\v{\delta}_i$ to reduce dimensions of the intra-regional OPF. 
The multi-area coordination problem \eqref{eq:maopf_theta} can be written in a compact form \eqref{eq:probP} as follows. 
\begin{subequations}
	\label{eq:probP}
	\begin{alignat}{2}
    &&  
	\underset{\v{x}_i, \v{\theta}_i}{\text{minimize}} & 
	\quad \sum_{i=1}^{N} \left ( \frac{1}{2} \v{x}_i^\T \v{H}_i \v{x}_i + \v{f}_i^\T \v{x}_i \right ), 
    \label{eq:probP.obj} \\ 
    && \text{subject to} &
	\quad \v{A}_i \v{x}_i \leq \v{b}_i + \v{C}_{ii} \v{\theta}_i + \sum_{j \in \nbd(i)} \v{C}_{ij} \v{\theta}_j, 
    \label{eq:probP.local} \\
	&&&
	\quad \v{\Theta}_i := \{ \v{D}_{ii} \v{\theta}_i + \sum_{j \in \nbd(i)} \v{D}_{ij} \v{\theta}_j \leq \v{r}_i \},
	\label{eq:probP.couple}	\\
    &&&
    \quad \forall  i = 1 , \dots, N 
	\notag
	\end{alignat}
\end{subequations}



\subsection{Preliminaries on CRE}
\label{subsec:prelim}

\begin{figure}[htbp]
    \centering
    \includegraphics[width=0.40\textwidth]{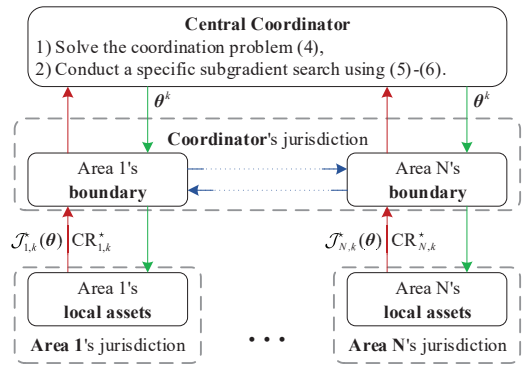} 
    \caption{The architecture of CRE. Red / green arrows indicate upward / downward communication links, and blue arrows are tie-lines among areas.}
    \label{fig:arch_cre}
\end{figure}

The compact form \eqref{eq:probP} can be interpreted as a primal decomposition formulation where $\v{x}_i, i = 1, \ldots, N$ are local variables, $\v{\theta}:= \left [ \v{\theta}_1^\T, \ldots, \v{\theta}_N^\T \right ]^\T$ are the coupling variables. 
The architecture of the primal decomposition method is shown in Figure \ref{fig:arch_cre}, where the local and coordination problems are iteratively solved.
Note that a central coordinator is needed, and it has jurisdiction over all areas' boundary phase angles, i.e., the coupling variables $\v{\theta}$. 

CRE also inherits the structure above. 
By leveraging the parametric programming theory, problem \eqref{eq:probP} can be formulated as a multi-parametric linear or quadratic programming (mp-LP/QP). The local problem $i$ parametric in $\v{\theta}$ is given as
\begin{subequations}
    \label{eq:mpP}
    \begin{alignat}{2}
    \Jcal_i^\star(\v{\theta}) := \ &
    \underset{\v{x}_i}{\text{minimize}} && 
    \ \frac{1}{2} \v{x}_i^\T \v{H}_i \v{x}_i + \v{f}_i^\T \v{x}_i, \label{eq:mpP.obj} \\ 
    & \text{subject to} &&
	\ \v{A}_i \v{x}_i \leq \v{b}_i + \v{C} \v{\theta}, \label{eq:mpP.cons}
    \end{alignat}
\end{subequations}
The properties of the value function $\Jcal^\star(\v{\theta}) = \textstyle \sum_{i=1}^N \Jcal^\star_i(\v{\theta})$ are summarized in the following lemma.
\begin{lemma}
    (see \cite{borrelli_bemporad_morari_2017}):
    \label{lemma.mpP}
    Consider the mp-LP/QP \eqref{eq:mpP}. The set of feasible parameters $\v{\Theta}^\star$ is a polyhedral set, i.e.,
    \begin{alignat*}{2}
        \v{\Theta}^\star := \cap_{i=1}^N \v{\Theta}^\star_i, \quad
        \v{\Theta}^\star_i := \{ \v{\theta} \ | \ \exists \v{x}_i: \v{A}_i \v{x}_i \leq \v{b}_i + \v{C} \v{\theta} \}
    \end{alignat*} 
    The value function $\Jcal^\star(\v{\theta}): \v{\Theta}^\star \rightarrow \Rset$ is convex piecewise linear / quadratic on a collection of polyhedral critical regions (CR), i.e., $\v{\Theta}^\star = \cap_{i=1}^N \cup_{k=1}^{K_i} \text{CR}^\star_{i,k}$, where $K_i$ is constant. 
\end{lemma}

Under Lemma \ref{lemma.mpP} and Figure \ref{fig:arch_cre}, CRE solves problem \eqref{eq:probP} by applying the following two steps recursively. 
\begin{enumerate}[leftmargin=*]
    \item \emph{Local evaluation:} Each area solves \eqref{eq:mpP} with $\v{\theta}^k$ to get $\Jcal^\star_{i,k}$, a slice of value function $\Jcal^\star_{i}$ defined over $\text{CR}^\star_{i,k}$.
    \item \emph{Coordination update:} Coordinator gathers all $\Jcal^\star_{i,k}$, $\text{CR}^\star_{i,k}$ to solve the problem
\end{enumerate}
    \begin{alignat}{2}
        \v{\theta}^\star := \{ \textstyle \argmin_{\v{\theta}} \sum_{i=1}^N \Jcal^\star_{i,k} \ | \ \cap_{i=1}^N \text{CR}^\star_{i,k}, \ \cap_{i=1}^N \v{\Theta}_i \}
        \label{eq:cre.coordination} 
    \end{alignat}
    and update $\v{\theta}^\star$ with a projected subgradient search 
\begin{alignat}{2}
    & \v{\theta}^{k+1} := \v{\theta}^\star - \varepsilon^{\text{stepsize}} \v{v}, \\
    & \v{v} := \{ \underset{\v{v}, \v{\eta}\geq \v{0}, \v{\zeta} \geq \v{0}}{\argmin} 
        \|\v{v}\|^2 \ | \  
        \v{v} = \partial \Jcal^\star \v{\eta} + \Ncal_{\v{\Theta}}^\star \v{\zeta}, 
        \ \bone^\T \v{\eta} = 1 \}
        \label{eq:CRE.LSQ}
\end{alignat}
Here, $\varepsilon^{\text{stepsize}}$ is constant, $\partial \Jcal^\star$, $\Ncal_{\v{\Theta}}^\star$ are the subdifferential and normal cone at $\v{\theta}^\star$ with $\v{\Theta} = \cap_{i=1}^N \v{\Theta}_i$, $\bone$ is an all one vector with proper dimensions. Lemma \ref{lemma.cre} shows the convergence of CRE.
\begin{lemma}
    (see \cite{guoRobustTieLineScheduling2018}):
    \label{lemma.cre}
    CRE obtains the optimal $\v{\theta}^\star$ after finite iterations when $\v{v}$ in \eqref{eq:CRE.LSQ} satisfies $\| \v{v} \| = 0$.
\end{lemma}



\section{Proposed RCDCRE Method}
\label{sec:rcdcre}

\subsection{Motivation and solution architecture}
\label{subsec:sol}

\begin{figure}[htbp]
    \centering
    \includegraphics[width=0.40\textwidth]{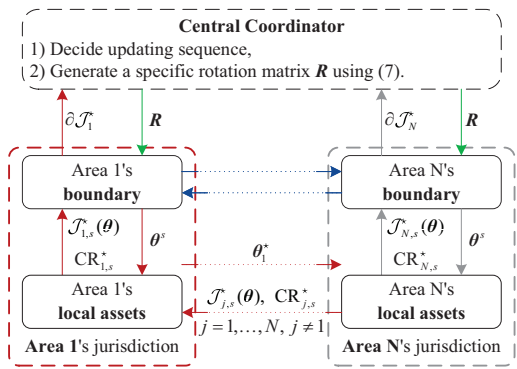} 
	\vspace{-7pt}
    \caption{The architecture of RCDCRE when area $i, i=1$ is updating.}
    \label{fig:arch_rcdcre}
\end{figure}

The coordinator in CRE gathers all local information and directly controls each area's boundary. 
In practice, system operators prefer weaker coordination and update boundary states asynchronously.
The proposed RCDCRE allows each area to reclaim jurisdiction of its boundary. As shown in Figure \ref{fig:arch_rcdcre}, the coordinator's jurisdiction is now intangible. Hence, we can deploy the coordination layer to any utility without privacy issues. 
Besides, recovering feasibility is a nontrivial task for all decomposition methods. 
The proposed RCDCRE method adopts penalty-based formulations to simultaneously take feasibility and optimality into account. Fast feasible recovery and convergence are achieved by implicitly exchanging the parametric mappings of the penalized intra-regional OPF.

Two main procedures of RCDCRE are highlighted in Figure \ref{fig:arch_rcdcre}, i.e., block coordinate descent (BCD) and coordinate system rotation (CSR). In RCDCRE, BCD allows each area to update boundary phase angles in a predefined sequence asynchronously (e.g., cyclic), shown in the areas' jurisdiction's dashed box. The current optimizing area (in the red dashed box) updates its boundary phase angles by locally interacting with other areas (in the grey dashed box). The information exchange includes the current area's boundary phase angles and parametric mappings of the present boundary state. All other areas' boundary phase angles are remained fixed. The internal dispatch of each region is independently adjusted according to the present boundary state. A stationary point $\theta^\star$ can be attained by the above BCD updates of the boundary phase angles. 
However, the convergence of $\theta^\star$ to optima, i.e., the global optimal boundary state of problem \eqref{eq:probP}, relies on two additional assumptions to the coordination problem \eqref{eq:cre.coordination}.
\begin{assumption}
    \label{assmp.BCD.obj}
    The objective function with respect to all block coordinates is differentiable, or it is not differentiable but the non-differentiable part is $N$ block-separable\cite{wrightCoordinateDescentAlgorithms2015}.
\end{assumption}

\begin{assumption}
    \label{assmp.BCD.cons}
    The feasible set is the direct sum of $N$ separable real Hilbert spaces of the block coordinates\cite{salzoParallelRandomBlockcoordinate2021}.
\end{assumption}

According to Lemma \ref{lemma.mpP}, for mp-LP/QP, $\Jcal^\star$ is a piecewise function, which is neither differentiable nor separable. The feasible set $\v{\Theta}^\star$ is a polyhedral set that is also non-separable. The BCD updates for boundary phase angles may get stuck at some non-optimal solutions.
Hence, an extra CSR step is introduced to remove convergence assumptions \ref{assmp.BCD.obj}-\ref{assmp.BCD.cons} when BCD updates the coordination problem \eqref{eq:cre.coordination}. 
For a stationary $\theta^\star$ after the BCD updates of all areas, the following rotation matrix $\v{R}$ is generated if $\theta^\star$ has not fulfilled the optimality condition. 
\begin{alignat}{2}
    \tilde{\v{v}} = \v{R} \v{v}, \quad \forall \ k:  \tilde{\v{v}} = \v{e}_k, \ 1 \leq k \leq \dim (\v{\theta})
    \label{eq:rcdcre.Givens}
\end{alignat}
where $\v{v}$ is a subgradient obtained from the solution to problem \eqref{eq:CRE.LSQ}. 
According to Lemma \ref{lemma.cre}, $\v{\theta}^\star$ is non-optimal when $ \|\v{v}\| \neq 0$.
The notation $\tilde{\v{v}}$ indicates any coordinatewise direction, where $\v{e}_k$ is the $k^{\textrm{th}}$ standard basis for the parameters $\v{\theta} \in \Rset^d$.
The matrix $\v{R}$ is what we want for conducting CSR, which can be obtained from compounding Givens matrices\cite[Section 5.1.8]{golubMatrixComputations2013}. 
The coordinator sends $\v{R}$ back to each area, the boundary phase angles' coordinate system can be rotated accordingly, as shown at the top of Figure \ref{fig:arch_rcdcre}.

The recovery of infeasible parameters comprises two steps. We first adopt a big$M$-penalty based formulation to each area $i$'s intra-regional OPF problem \eqref{eq:mpP}.
\begin{subequations}
    \label{eq:mpP_bigM}
    \begin{alignat}{2}
    \Jcal_i^\star(\v{\theta}) := \ &
    \underset{\v{x}_i}{\text{minimize}} && 
    \ \bar{\v{x}}_i^\T \bar{\v{H}}_i \bar{\v{x}}_i + \bar{\v{f}}_i^\T \bar{\v{x}}_i, \\ 
    & \text{subject to} &&
	\ \bar{\v{A}}_i \bar{\v{x}}_i \leq \bar{\v{b}}_i + \bar{\v{C}}_i \v{\theta},
    \end{alignat}
\end{subequations}
where the variables and coeficients are given by
\begin{alignat*}{2}
	\bar{\v{x}}_i = \begin{bmatrix}
        \v{x}_i \\
        \v{s}_i
    \end{bmatrix}, \quad  
    \bar{\v{H}}_i = \begin{bmatrix}
        \v{H}_i & \v{0} \\
        \v{0} & \v{0}
    \end{bmatrix}, \quad  
    \bar{\v{f}}_i = \begin{bmatrix}
        \v{f}_i \\
        M \bone \\
    \end{bmatrix}, \\
    \bar{\v{A}}_i = \begin{bmatrix}
        \v{A}_i & -\v{I} \\
        \v{0}  & -\v{I} 
    \end{bmatrix}, \quad
    \bar{\v{b}}_i = \begin{bmatrix}
        \v{b}_i \\
        \v{0}
    \end{bmatrix}, \quad
    \bar{\v{C}}_i = \begin{bmatrix}
        \v{C}_i \\
        \v{0}
    \end{bmatrix}.
\end{alignat*}
The introduction of slack variable $\v{s}_i$ ensures an always feasible solution to problem \eqref{eq:mpP_bigM} for any given parameters $\v{\theta}$. 
Note that to reduce the number of $\v{s}_i$, it is recommended to remove the redundant constraints to the jointly feasible region of $\v{x}_i$-$\v{\theta}$ space before the reformulation.

To adjust the tie-line schedules that physically connected to area $i$, area $i$ solve an interregional OPF problem \eqref{eq:rcdcre_mp_L1Reg} to update its boundary phase angles. Problem \eqref{eq:rcdcre_mp_L1Reg} can be viewed as problem \eqref{eq:cre.coordination} with $\ell_1$-penalization to the $\cap_{i=1}^N \v{\Theta}_i$ and fixed block coordinates to other area's boundary state.
\begin{subequations}
	\label{eq:rcdcre_mp_L1Reg}
	\begin{alignat}{2}
    &&  
    \underset{\v{\theta}}{\text{minimize}} & \quad 
    \sum_{i=1}^{N} \Jcal^\star_{i,s}
    + \sigma \bone^T \max \left \{ \v{D}_{i} \v{\theta} - \v{r}_{i}, \v{0} \right \}
	\label{eq:rcdcre_mp_L1Reg.obj} \\
	&& \text{subject to} &
    \quad \cap_{i=1}^N \text{CR}^\star_{i,s}, 
    \quad \v{\theta}_{\lnot i} = \v{\theta}^\star_{\lnot i}
	\label{eq:rcdcre_mp_L1Reg.CR} 
	\end{alignat}
\end{subequations}
Here, $\v{D}_{i}$, $\v{r}_{i}$ denote compact coefficients of $\v{\Theta}_i$, $\sigma$ is a penalty coefficient. The boundary phase angles that are not controlled by area $i$ are denoted by $\v{\theta}_{\lnot i}$.
By Lemma \ref{lemma.mpP}, the value function segment $\Jcal^\star_{i,s}$ and critical region $\text{CR}^\star_{i,s}$ can be explicitly expressed as
\begin{subequations}
	\label{eq:mpP_VF_CR}
	\begin{alignat}{2}
    &&  
    \Jcal^\star_{i,s} & :=  
    \frac{1}{2} \v{\theta}^\T \hat{\v{H}}_{i,s} \v{\theta} + \hat{\v{f}}_{i,s}^\T \v{\theta} + \hat{c}_{i,s}, \\
    &&  
    \text{CR}^\star_{i,s} & := \{ \v{\theta} \ | \ 
        \hat{\v{D}}_{i, s} \v{\theta} \leq \hat{\v{r}}_{i,s}
        \}.     
	\end{alignat}
\end{subequations}
Under non-degenerate assumptions, \eqref{eq:mpP_VF_CR} can be directly obtained from KKT conditions. See \cite{guoRobustTieLineScheduling2018, guoCoordinatedMultiAreaEconomic2017} for details. 
When the primal or dual solutions are non-unique, \eqref{eq:mpP_VF_CR} can still be obtained with additional techniques, such as (a) projecting from a higher dimensional solution space \cite{jonesPolyhedralProjectionParametric2008,akbariImprovedMultiparametricProgramming2018a}, (b) solving auxiliary problems of the optimal set \cite{jonesLexicographicPerturbationMultiparametric2007,spjotvoldContinuousSelectionUnique2007}, or (c) searching the combinational tree of active constraints \cite{guptaNovelApproachMultiparametric2011a}.

The $\ell_1$-penalty can return a feasible solution outside the initial feasible region. Hence, it allows RCDCRE to take an arbitrary initial $\v{\theta}$  as a cold start. 
By exchanging parametric mappings of the present boundary state, each area updates its boundary phase angles locally in a BCD fashion by solving penalty-based intra-regional OPF \eqref{eq:mpP_bigM} and interregional OPF \eqref{eq:rcdcre_mp_L1Reg}. The central coordinator conducts CSR if all areas finish BCD updates and the optimal condition has not been met.
The proposed RCDCRE converges to optimal by iteratively applying the above steps.

Compared to existing techniques, LR decomposition needs feasible recovery steps for general problems\cite{gustavssonPrimalConvergenceDual2015}. Large penalty factor decreases ADMM's performance, yet the intermediate solution is still asymptotic feasible. 
Benders decomposition with feasible cuts cannot track varied system conditions, and the convergence is poor if it adopts a big$M$-penalty formulation \cite{candasComparativeStudyBenders2020}. 
Our method considers an exact parametric relation between penalty formulations and feasibility.
Section \ref{sec:case} shows that it can recover feasible solutions in a few iterations from infeasible boundary phase angles.


Algorithm \ref{alg:RCDCRE} shows the details of RCDCRE. The algorithm sets an initial value $\theta^\star$. By defaut, $\theta^\star$ is not the coordinatewise optimal point to all the areas, i.e., $\Aset_w \gets \{ 1, \ldots, N\}$. List $\Aset_o$ is the complement of $\Aset_w$. The coordinates in area $i$ are saved in $\Eset_{i, w}$. For example, suppose area $i$ controls $p^{\textrm{th}}$ and $q^{\textrm{th}}$ of the boundary phase angle $\v{\theta}$, then we have
\begin{alignat}{2}
	\Eset_{i, w} = \{ e_p, \ -e_p, \  e_q,  \  -e_p \} \label{eq:BCD.coordinates}
\end{alignat}
where, $e_p$ is the $p^{\textrm{th}}$ basis, $-e_p$ is its negative direction. Set $\Eset_{i, o}$ is the complement of $\Eset_{i, w}$ for area $i$.

Step \ref{alg:RCDCRE.BCD.ini} obtain the current optimizing area and boundary state.
Steps \ref{alg:RCDCRE.CRCRE.judge}-\ref{alg:RCDCRE.BCD.index_update3} update the local boundary phase angles. 
The coordinatewise search is conducted in step \ref{alg:RCDCRE.BCD.search}. The new explored parameter is denoted as $\v{\theta}^s$. 
In steps \ref{alg:RCDCRE.BCD.map1}-\ref{alg:RCDCRE.BCD.map2}, all areas solve the big$M$-penalty intra-regional OPF problem \eqref{eq:mpP_bigM} to dispatch internal assets given the explored parameter $\v{\theta}^s$.
Area $i$ then gathers the parametric mappings of the present boundary state. Step \ref{alg:RCDCRE.BCD.coordination} solves the $\ell_1$-penalty interregional OPF problem with other areas' boundary phase angles $\v{\theta}_{\lnot i}$ are fixed. The current optimal solution is denoted as $\v{\theta}^\textrm{tmp}$.

The subdifferential $\partial \Jcal^\star_i$, normal cone $\Ncal_{\v{\Theta}_i}^\star$, and explored coordinate $s$ are updated in steps \ref{alg:RCDCRE.BCD.index_update1}-\ref{alg:RCDCRE.BCD.index_update3} when area $i$ optimizes the $\v{\theta}_i$.\footnote{Note $\Ncal_{\v{\Theta}_i}^\star$ has the same updating rule as $\partial \Jcal^\star_i$, and is omitted for brevity.}
If area $i$ has explored all the coordinate directions, then switch to another area in $\Aset_w$.
The coordinator gathers all $\partial \Jcal^\star_i$, $\Ncal_{\v{\Theta}_i}^\star$ for judging optimality, and updates $\Aset_w$ for conducting BCD, shown in steps \ref{alg:RCDCRE.subgrad.judge}-\ref{alg:RCDCRE.subgrad.append}.
Step \ref{alg:RCDCRE.check.BCD} checks whether all area finishes the BCD update for their boundary phase angles. If true, then steps \ref{alg:RCDCRE.LSQ}-\ref{alg:RCDCRE.return} do the global optimality check. 
For the non-optimal $\theta^\star$, CSR is conducted in steps \ref{alg:RCDCRE.rotation.generation}-\ref{alg:RCDCRE.rotation.reform2}. Otherwise, Algorithm terminates with the optimal global solutions.

\begin{algorithm}[!htbp]
    \caption{Penalty based RCDCRE method.}
    \label{alg:RCDCRE}

    Set 
    $\v{\theta}^\star \gets$ a start point,
    $\Aset_w \gets \{ 1, \ldots, N\}$, 
    $\Aset_o \gets \emptyset$.

    \For{$k=1, 2, \ldots$}{
    
    $i \gets \Aset_w(1)$, $\Aset_w \gets \Aset_w \backslash i$, $\Aset_o \gets \Aset_o \cup i$,
    $\v{\theta}^\star_{i} \gets \v{\theta}^\star$, 
    $\Eset_{i, w} \gets$ \eqref{eq:BCD.coordinates}, $\Eset_{i, o} \gets \emptyset$
    \label{alg:RCDCRE.BCD.ini}

    \While(){$\Eset_{i, w} \neq \emptyset$}{    
        \label{alg:RCDCRE.CRCRE.judge}

    $\v{v} \gets \Eset_{i, w}(1)$, 
    $\v{\theta}^{s} \gets \v{\theta}^\star_{i} - 
        \varepsilon^{\text{stepsize}} \v{v}$
    \label{alg:RCDCRE.BCD.search}
        
    \For{$j=1, \ldots, N$}{
        \label{alg:RCDCRE.BCD.map1}

        $\v{x}^\star_j \gets $ solution to problem \eqref{eq:mpP_bigM} with $\v{\theta}^{s}$
        
        $\Jcal^\star_{j,s}$, $\text{CR}^\star_{j,s} \gets $ Parametric map of \eqref{eq:mpP_VF_CR}
        \label{alg:RCDCRE.BCD.map2} 
    }

    $\v{\theta}^\textrm{tmp} \gets $ solution to problem \eqref{eq:rcdcre_mp_L1Reg}
    \label{alg:RCDCRE.BCD.coordination}

	\eIf{$ \| \v{\theta}^\textrm{tmp} - \v{\theta}^\star_{i} \| \geq \varepsilon^{\text{stepsize}}$}{
        \label{alg:RCDCRE.BCD.index_update1}
        $\Eset_{i, w} \gets \Eset_{i, w} \cup \Eset_{i, o}$, $\Eset_{i, o} \gets \emptyset$,
        $\v{\theta}^\star_{i} \gets \v{\theta}^\textrm{tmp}$, 
        $\partial \Jcal^\star_i \gets \nabla \Jcal^\star$
        \label{alg:RCDCRE.BCD.index_update2}
    }{
        $\Eset_{i, w} \gets \Eset_{i, w} \backslash \v{v}$, $\Eset_{i, o} \gets \Eset_{i, o} \cup \v{v}$,
        $\partial \Jcal^\star_i \gets \partial \Jcal^\star_i \cup \nabla \Jcal^\star$
        \label{alg:RCDCRE.BCD.index_update3}		
	}    
    
    }

	\eIf{$ \| \v{\theta}^\star - \v{\theta}^\star_{i} \| \geq \varepsilon^{\text{stepsize}}$}{
        \label{alg:RCDCRE.subgrad.judge}

        $\v{\theta}^\star \gets \v{\theta}^\star_{i}$, 
		$\partial \Jcal^\star \gets \{ \partial \Jcal^\star_i \}$,
        $\Aset_w \gets \Aset_w \cup \Aset_o$
        \label{alg:RCDCRE.BCD.new}

	}{
		$\partial \Jcal^\star \gets \partial \Jcal^\star \cup  \{ \partial \Jcal^\star_i \}$
        \label{alg:RCDCRE.subgrad.append}
	}

    \If(){$\Aset_w = \emptyset$}{
    \label{alg:RCDCRE.check.BCD}

	$\v{v} \gets $ solution to the subgradient direction \eqref{eq:CRE.LSQ}
    \label{alg:RCDCRE.LSQ}

	\eIf{$ \|\v{v}\| \geq \varepsilon^{\text{optimal}} $}{
        \label{alg:RCDCRE.check.opt}

        $\v{R} \gets$ rotation matrix from \eqref{eq:rcdcre.Givens} 
        \label{alg:RCDCRE.rotation.generation}

        \For{$i=1, \ldots, N$}{
            \label{alg:RCDCRE.rotation.reform1}
            Reform problem \eqref{eq:mpP_bigM} using $\v{R}$
            \label{alg:RCDCRE.rotation.reform2}
    
        }

	}{

        \Return $\v{\theta}^\star$
        \label{alg:RCDCRE.return}
	}
    }

    }
\end{algorithm}

Theorem \ref{thm.RCDCRE} states the convergence of Algorithm \ref{alg:RCDCRE}.
\begin{theorem}
    \label{thm.RCDCRE}
    RCDCRE converges to optimal parameters after finite BCDs and CSRs.
\end{theorem}

\subsection{Illustrative example}
\begin{figure}[htbp]
    \centering   
    \includegraphics[width=0.45\textwidth]{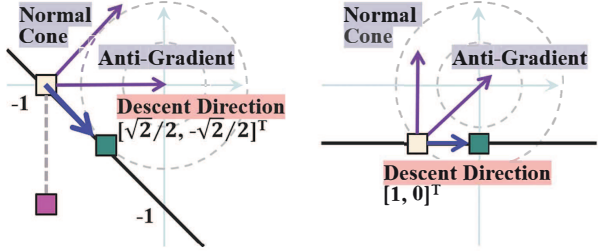}  
	\caption{An example of coordinate system before and after rotation.}
    \label{fig:rcdcre_toy}
\end{figure}

We show how CSR avoids the stucking of BCD. 
Consider problem \eqref{eq:Toy} with two agents and an initial point $(-1, -1)$.
\begin{subequations}
    \label{eq:Toy}
    \begin{alignat}{2}
    &
    \underset{\theta_1,\theta_2}{\text{minimize}} && 
    \quad \theta_1^2 + \theta_2^2, 
    \label{eq:Toy.obj} \\ 
    &
    \text{subject to} && \quad \theta_1 + \theta_2 \leq -1
    \label{eq:Toy.cons}
    \end{alignat}
\end{subequations}
In the left of Figure \ref{fig:rcdcre_toy}, $\theta_2$ do the coordinatewise minimization first, the optimal solution is $(-1, 0)$. 
However, $\theta_1$ cannot find any further improvement from $(-1, 0)$ to the global optima $(-0.5, -0.5)$, as \eqref{eq:Toy.cons} violates Assumption \ref{assmp.BCD.cons}.

Notice at $(-1, 0)$, a descent direction $\left [\sqrt{2}/2,\ -\sqrt{2}/2 \right ]^\T$ can be induced by the set plus of normal cone and gradient. 
Hence, we can use \eqref{eq:rcdcre.Givens} to rotate it into $\left [1,\ 0 \right ]^\T$, corresponding to the standard basis of $\theta_1$.
The transformation is given by
\begin{alignat}{2}
    \begin{bmatrix}
        1 \\
        0 \\
    \end{bmatrix} = 
    \begin{bmatrix}
        \frac{\sqrt{2}}{2} & -\frac{\sqrt{2}}{2}\\
        \frac{\sqrt{2}}{2} & \frac{\sqrt{2}}{2} \\
    \end{bmatrix}
    \begin{bmatrix}
        \frac{\sqrt{2}}{2} \\
        -\frac{\sqrt{2}}{2} \\
    \end{bmatrix}
    \label{eq:Toy.rot}
\end{alignat}
where the matrix in \eqref{eq:Toy.rot} is the rotation matrix $\v{R}$. 
By substituting $\v{\theta} = \v{R}^{-1} \tilde{\v{\theta}}$ into problem \eqref{eq:Toy}, the original problem is equivalently reformed in the new coordinate system as 
\begin{alignat*}{2}
    \underset{\tilde{\theta}_1,\tilde{\theta}_2}{\text{minimize}}
    \quad \tilde{\theta}_1^2 + \tilde{\theta}_2^2, \quad
    \text{subject to} \quad \tilde{\theta}_2 \leq -\frac{\sqrt{2}}{2}
\end{alignat*}
As shown in the right of Figure \ref{fig:rcdcre_toy}, current descent direction becomes $\left [1,\ 0 \right ]^\T$, which can be solved by BCD readily.


\section{Case Studies}
\label{sec:case}

We used MATLAB 2021a with CPLEX version 12.9.0 to conduct all the simulations. 
Network data were obtained from the MATPOWER 7.1 \cite{Zimmerman2011MATPOWERSO}. 
For the illustrative purpose, feasible dispatch region $\v{\Theta} = \cap_{i=1}^2 \v{\Theta}_i$ in the two-area 44-bus system was divided into critical regions using MPT 3.0 \cite{MPT3}.
A four-area 472-bus network was designed for comparative analysis. 
Unless specified, a zero cold start is adopted for the initial boundary phase angles $\v{\theta}$.

\subsection{Simulation on a two-area 44-bus network}

\begin{figure}[htbp]
    \centering
    \includegraphics[width=0.4\textwidth]{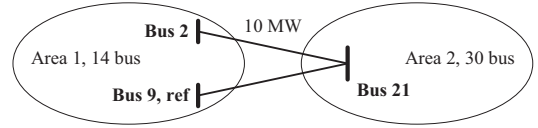} 
    \caption{Two-area 44-bus network.}
    \label{fig:2A44B2TL.0}
\end{figure}

\begin{figure*}[htbp]
    \centering
    \subfloat[]{ 
        \includegraphics[width=0.32\textwidth]{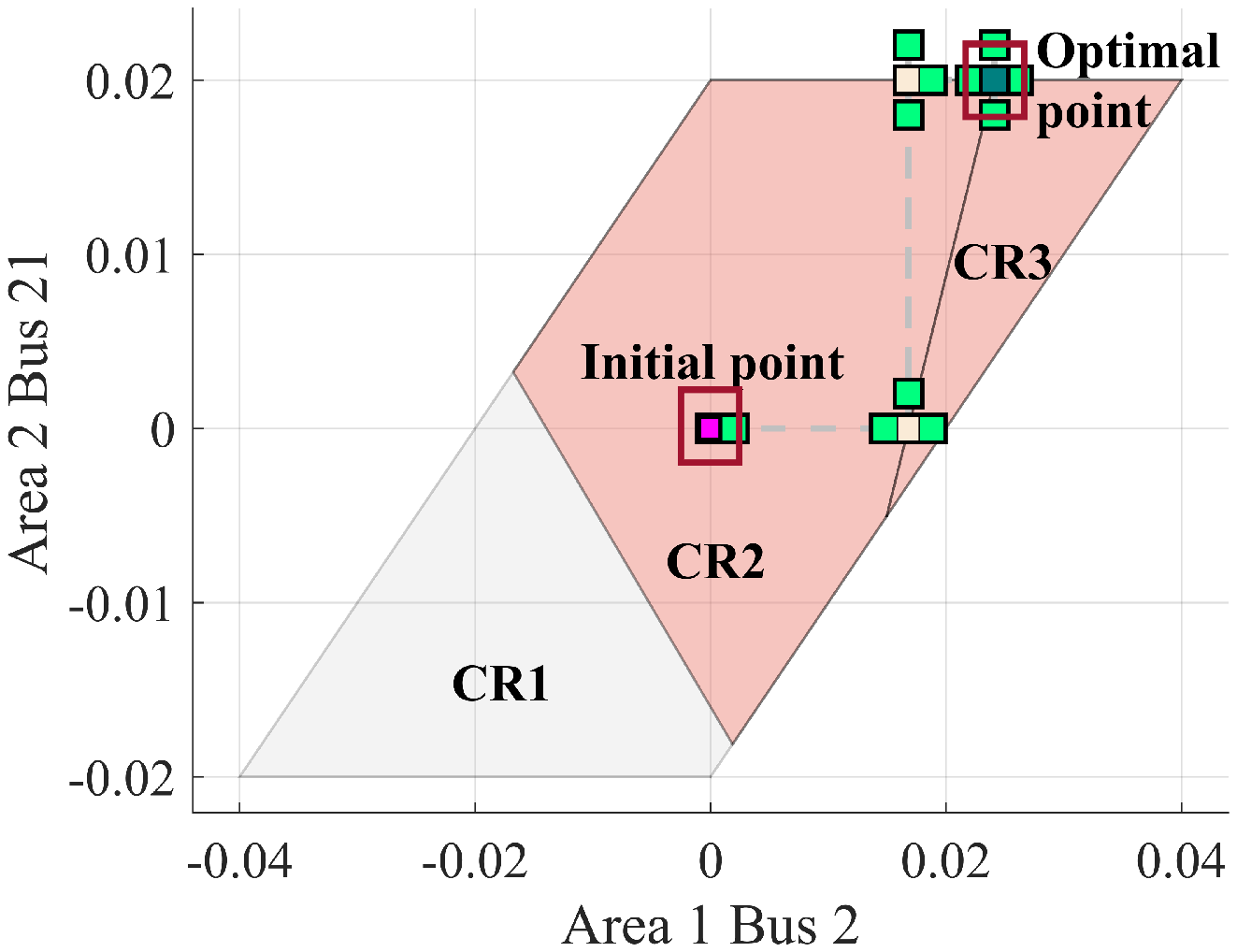} 
		\label{fig:2A44B2TL.1}
    }
    \subfloat[]{ 
        \includegraphics[width=0.32\textwidth]{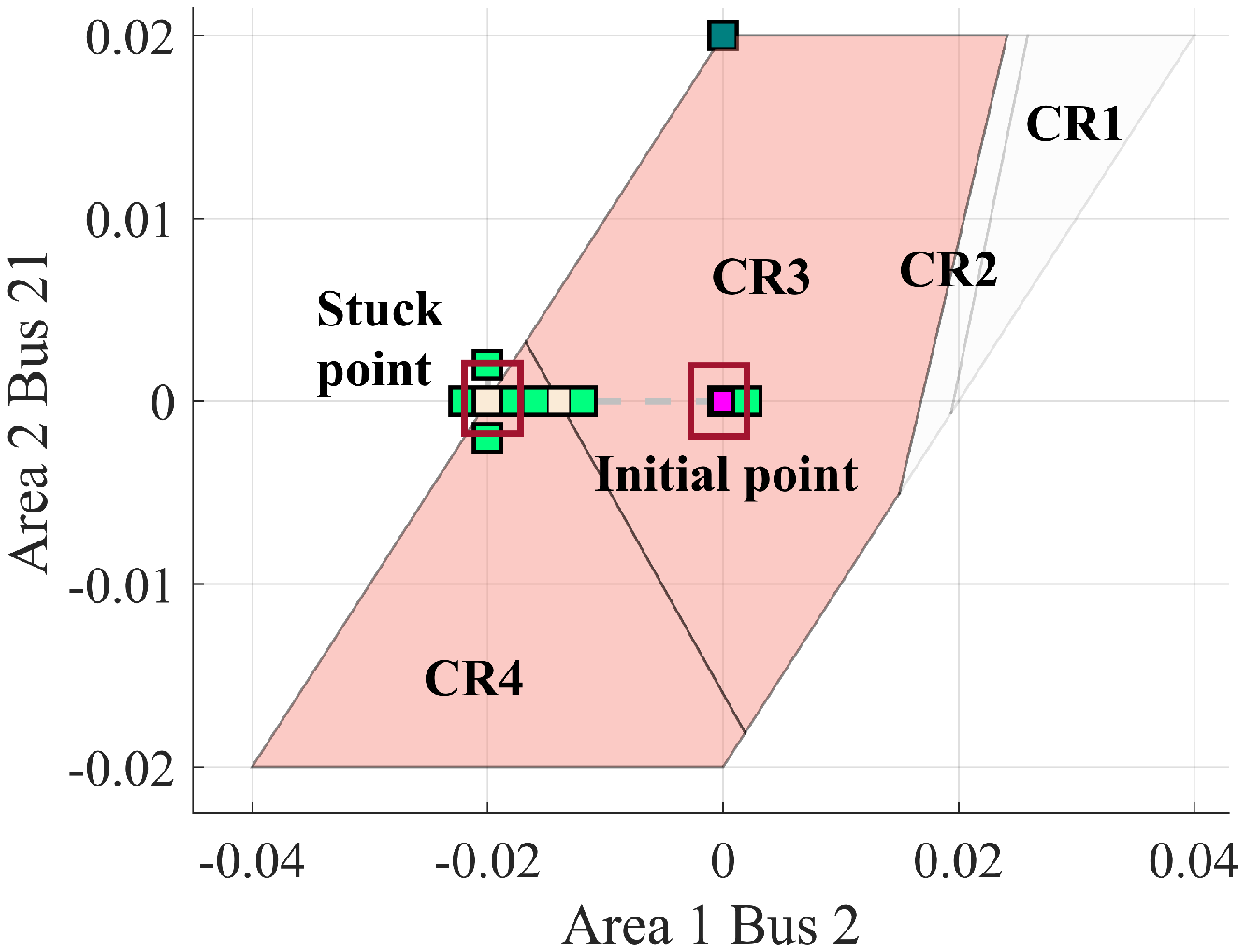} 
		\label{fig:2A44B2TL.2}
    }
    \subfloat[]{ 
        \includegraphics[width=0.32\textwidth]{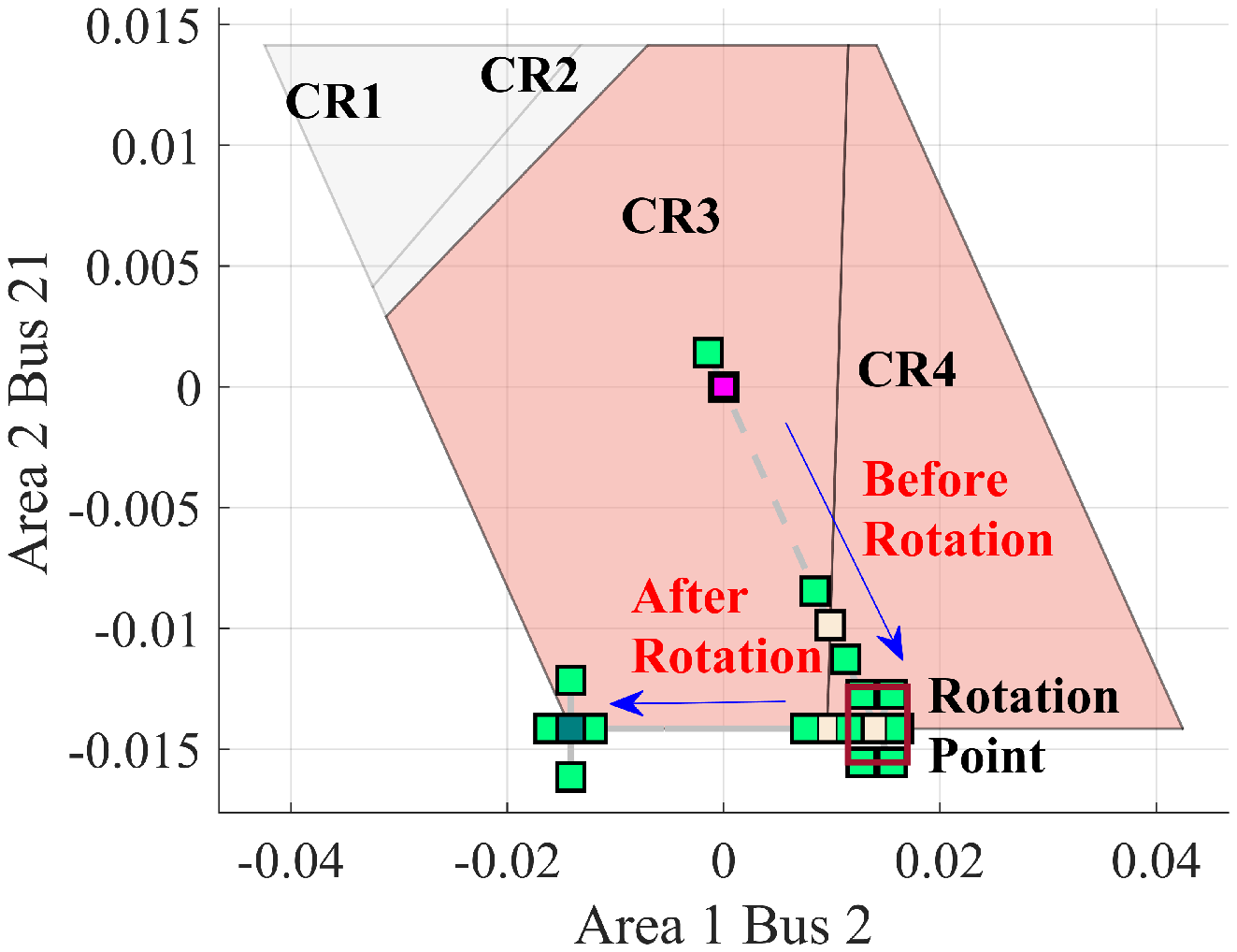} 
		\label{fig:2A44B2TL.3}
    }    
	\caption{Simulations on two-area 44-bus network (a) Linear generation costs, (b)-(c) Quadratic generation costs before and after rotation.}
\end{figure*}

Consider the two-area power system in Figure \ref{fig:2A44B2TL.0} that stitched IEEE 14 and 30 systems together with 2 tie-lines as shown. All tie-line capacities were set as $10$ MW, and internal line capacities were $100$ MW. Linear cost coefficients $\v{c}_i$ of the generators were perturbed to $\widetilde{\v{c}}_i := \v{c}_i \circ \left( 0.99 + 0.02 \v{\xi}_i \right)$, for $i = 1,2$, where entries of $\v{\xi}_i$ are independent $\Ncal(0,1)$ (standard normal) variables. 
Simulation results that neglect quadratic generation costs are shown in Figure \ref{fig:2A44B2TL.1}. The proposed RCDCRE method can find the exact optimal solutions after switching the coordinates four times. The feasible space $\v{\Theta}$ contains three critical regions, and only two of them need to be explored for the entire RCDCRE procedure. If we do not neglect the quadratic costs, then $\v{\Theta}$ comprises four critical regions. As shown in Figure \ref{fig:2A44B2TL.2}, the problem stuck at a non-optimal point after switching coordinates two times. After the CSR, RCDCRE can find the optima by another two switchings, shown in Figure \ref{fig:2A44B2TL.3}. 

\begin{figure}[htbp]
    \centering
    \includegraphics[width=0.48\textwidth]{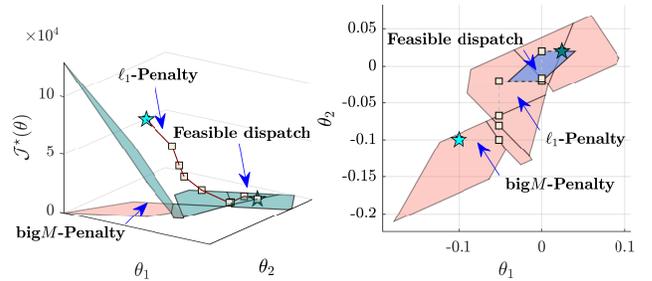}         
    \vspace{-7pt}
	\caption{Simulations with linear costs and start point $(-0.1, -0.1)$.}
    \label{fig:2A44B2TL.inf.1}
\end{figure}

By changing the start point to $(-0.1, -0.1)$, the trajectory of feasibility recovery is illustrated in Figure \ref{fig:2A44B2TL.inf.1}. Here, area 1 bus 2 is denoted by $\theta_1$, and $\theta_2$ represents area 2 bus 21. We adopt a linear generation cost for illustration. 
Note feasible parameter space $\v{\Theta}^\star = \Rset^2$, as any parameter $\v{\theta}$ has at least a feasible solution to \eqref{eq:mpP_bigM}. 
Hence, by lemma \ref{lemma.cre}, the value function $\Jcal^\star$ is a collection of piecewise affine hyperplanes in $\Rset^2$, as shown in the left of Figure \ref{fig:2A44B2TL.inf.1}. 
The converging process roughly comprises three stages for this case. In the first stage, the intra-regional OPF is infeasible to \eqref{eq:mpP} but feasible to \eqref{eq:mpP_bigM}. Hence, RCDCRE returns a large objective defined over that critical region, mainly caused by big$M$-penalty. After one BCD update along $\Jcal^\star$, boundary phase angles $\v{\theta}$ become feasible to \eqref{eq:mpP} for all areas rapidly.
The $\ell_1$-penalty to \eqref{eq:rcdcre_mp_L1Reg} maximizes its effect in the second stage. We can observe that $\ell_1$-penalty allows several critical regions generated outside $\v{\Theta}$, shown in the shaded area in the right of Figure \ref{fig:2A44B2TL.inf.1}. Within switching coordinate two times, $\v{\theta}$ converge back to $\v{\Theta}$. The final stage, when $\v{\theta} \in \v{\Theta}$ and feasible to \eqref{eq:mpP} is identical to the trajectory in Figure \ref{fig:2A44B2TL.1}. 
Similar convergence also applies to the quadratic generation cost case. 

\subsection{Simulation on a four-area 472-bus network}

\begin{figure}[htbp]
    \centering
    \includegraphics[width=0.45\textwidth]{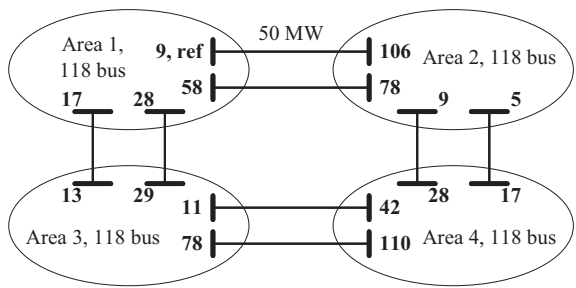} 
	\caption{Four-area 472-bus network.}   
    \label{fig:4A472B4TL.0} 
\end{figure}

We next show the results on a large four-area 472-bus network with $50$ MW tie-line and $500$ MW internal line capacities. The topology is depicted in Figure \ref{fig:4A472B4TL.0}. A quadratic generation cost was adopted, where the linear coefficients were perturbed similar to the two-area one.

\begin{figure}[htbp]
    \centering
    \includegraphics[width=0.48\textwidth]{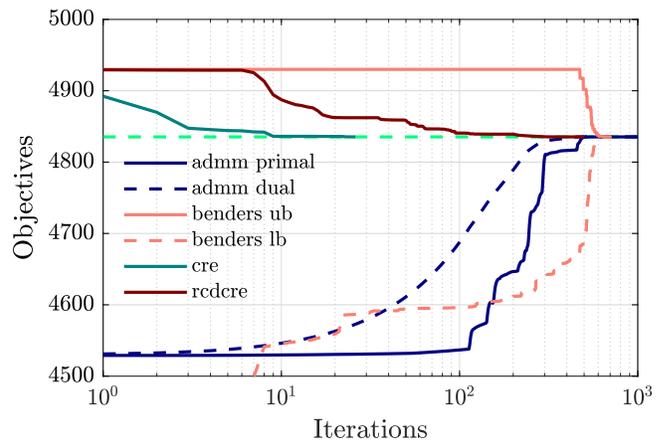} 
	\caption{Comparison in objective convergence of four algorithms.} 
    \label{fig:4A472B4TL.1}   
\end{figure}

RCDCRE is compared with three other approaches, i.e., CRE taken from \cite{guoRobustTieLineScheduling2018}, Consensus ADMM taken from \cite[Algorithm 2]{kargarianDistributedDecentralizedDC2018}, and Generalized Benders decomposition modified from \cite[Section 5.1]{birgeIntroductionStochasticProgramming2011}.
The objective convergence results are shown in Figure \ref{fig:4A472B4TL.1}. 
Compared to the centralized solution, the CRE method converges to the exact solution after 26 iterations. RCDCRE takes 114 iterations converging to the $10^{-3}$ relative optimality gap. Meanwhile, it needs 477 and 589 iterations for ADMM primal objective and Benders upper bound converging to the same accuracy level. RCDCRE allows each area to update its boundary states rather than the coordinated update like CRE and Benders. Moreover, it enjoys a better convergence over ADMM and Benders decomposition due to adequate information exchanged.

\section{Conclusions}
\label{sec:Conclusions}

We have proposed RCDCRE as a distributed privacy-preserving and cold start approach to solving the multi-area OPF for the first time. 
RCDCRE does not need a coordinator to control and update the boundary phase angles synchronously as in the Benders and CRE method. 
Parametric penalty-based formulations are adopted to achieve rapid feasibility recovery and convergence.
The limitation of the RCDCRE method is that we still require coordination to do the coordinate system rotation to ensure convergence. 
Nevertheless, RCDCRE shows finite convergent property on a two-area 44-bus network. 
Results on a four-area 472-bus network show that RCDCRE enjoys a faster convergence rate than the ADMM and Benders decomposition methods.
The structure and simulations results indicate that RCDCRE is favorable to the operation of interconnected power systems, and it also applies to generic multi-agent optimization problems. 
Future work will entail degeneracy handling, coordination reduction, and applying the proposed method to general convex problems in power system operations.

\bibliographystyle{IEEEtran}
\bibliography{IEEEabrv,ref_cdc}

\appendix


\subsection{Proof of Theorem \ref{thm.RCDCRE}}
\label{secsec:proofs}

To prove theorem \ref{thm.RCDCRE}, we first impose an assumption about the stepsize in step \ref{alg:RCDCRE.BCD.search} of algorithm \ref{alg:RCDCRE}.
\begin{assumption}
    \label{assmp.stepsize}
    Given $\varepsilon^{\text{stepsize}}$, $\v{\theta}^s$ can identify but not step across an adjacent critical region if $\v{\theta}^\star$ is on the boundary of current critical region $\cap_{i=1}^N \text{CR}^\star_{i,s}$. Otherwise, $\varepsilon^{\text{stepsize}}$ ensures that $\v{\theta}^s$ does not step out of $\cap_{i=1}^N \text{CR}^\star_{i,s}$.
\end{assumption}

If $\| \v{\theta}^\textrm{tmp} - \v{\theta}^\star_{i} \| \geq \varepsilon^{\text{stepsize}}$, such a stepsize can ensure $\Jcal( \v{\theta}^\textrm{tmp} ) < \Jcal(\v{\theta}^\star_{i})$. Hence, the $\Jcal(\v{\theta}^\star)$ during the solution process is always non-increasing.

We proof Theorem \ref{thm.RCDCRE} in three steps.
\begin{enumerate}[leftmargin=*]
    \item Without penalty formulation, RCDCRE converges when first-order non-differentiable optimality condition is fulfilled for feasible parameters $\v{\theta} \in \v{\v{\Theta}}^\star$. \label{proof:CRE}
    \item The big$M$-penalty formulation can return the same optimal solution when initial $\v{\theta}$ is infeasible to \eqref{eq:mpP}. \label{proof:bigM}
    \item The $\ell_1$-penalty formulation returns the same optimal solutions when starting from $\v{\theta} \notin \v{\v{\Theta}}$. \label{proof:L1}
\end{enumerate} 

\subsubsection{Proof of \ref{proof:CRE})}

For block coordinate descent (BCD), it is well known that it converges to optimal in succesive blockwise minimizations when assumptions \ref{assmp.BCD.obj}-\ref{assmp.BCD.cons} holds. 
However, when Assumptions \ref{assmp.BCD.obj}-\ref{assmp.BCD.cons} fails, the CSR conducts in RCDCRE is given by
\begin{subequations}
    \label{eq:rcdcre_rot}
    \begin{alignat}{2}
        & \tilde{\v{v}} = \v{R} \v{v}, \quad \tilde{\v{\theta}}^\star = \v{R} \v{\theta}^\star, \\
        & \tilde{\v{C}}_{i} = \v{C}_{i} \v{R}^\T, \quad
        \tilde{\v{D}}_{i} = \v{D}_{i} \v{R}^\T, \quad \forall i = 1, \ldots, N. 
    \end{alignat}
\end{subequations}
where $\v{R}$ is the rotation matrix, the rotated parameters and coefficients are denoted in $\left [^\sim \right ]$ values. Note $\v{R}^\T = \v{R}^{-1}$ as $\v{R}$ is orthogonal. The subgradinet direction $\v{v}$ is obtained from \eqref{eq:CRE.LSQ}. The rotated subgradient direction is denoted by $\tilde{\v{v}}$, 
and it is also a coordinate direction according to $\v{R}$. 
Hence, an exploration of $\tilde{\v{\theta}}^\star$ in step \ref{alg:RCDCRE.BCD.search} of Algorithm \ref{alg:RCDCRE} is given by
$$\v{\theta}^s = \tilde{\v{\theta}}^\star - \varepsilon^{\text{stepsize}} \tilde{\v{v}}$$

\begin{remark}
    \label{rmk:rot}
    The rotated subgradient direction $\tilde{\v{v}}$ is also the subgradient in the rotated coordinate system.
\end{remark}

\begin{proof}
    According to \eqref{eq:CRE.LSQ}, remark \ref{rmk:rot} implies the explored subdifferential and normal cones are also the subdifferential and normal cones in the rotated system. For each critical region with active and inactive set partitions $\Acal$, $\Ical$, by KKT conditions we have
    \begin{subequations}
        \label{eq:mpP_KKT}
        \begin{alignat}{2}
            &&& \v{H} \v{x} + \v{f} + 
            \v{A}_\Acal^\T \v{\lambda}_\Acal = \v{0}, \\
            &&& \v{A}_\Acal \v{x} = \v{b}_\Acal + \tilde{\v{C}}_\Acal \tilde{\v{\theta}}.
        \end{alignat}
    \end{subequations}     
    WLOG, we assume $\v{H} = 0$, the problem reduce to a mp-LP, and subscript $i$ has been dropped for simplicity. The $\v{C}_\Acal \v{\theta}$ is replaced by $\tilde{\v{C}}_\Acal \tilde{\v{\theta}}$ because of \eqref{eq:rcdcre_rot}. If $\v{A}_\Acal$ is full rank, then we have the primal mapping and value function as follows
    \begin{subequations}
        \begin{alignat}{2}
            \v{x}^\star(\tilde{\v{\theta}}) = \v{A}_\Acal^{-1} (\v{b}_\Acal + \tilde{\v{C}}_\Acal \tilde{\v{\theta}}), \\
            \Jcal^\star(\tilde{\v{\theta}}) = \v{f}^\T \v{A}_\Acal^{-1} (\v{b}_\Acal + \tilde{\v{C}}_\Acal \tilde{\v{\theta}}).
        \end{alignat}
    \end{subequations}     
   The gradient of current value function segment is
    \begin{subequations}
        \begin{alignat}{2}
            && \nabla \Jcal^\star(\tilde{\v{\theta}}) 
            & = \tilde{\v{C}}_\Acal^\T (\v{A}_\Acal^{-1})^\T \v{f}, \\
            &&& = \v{R} \v{C}_\Acal^\T (\v{A}_\Acal^{-1})^\T \v{f}, \\
            &&& = \v{R} \nabla \Jcal^\star(\v{\theta}).
        \end{alignat}
    \end{subequations}
    Similar argument also applies for mp-QP, degenerate conditions, as well as the rotated normal cones. By substituting $\partial \Jcal^\star(\tilde{\v{\theta}})$, $\Ncal^\star_{\tilde{\v{\Theta}}}(\tilde{\v{\theta}})$ into \eqref{eq:CRE.LSQ}, we can obtain $\tilde{\v{v}}$ is the optimal solution to \eqref{eq:CRE.LSQ}, which completes the proof.
\end{proof}

Under remark \ref{rmk:rot}, the coordinatewise search is also a subgradient search, and we have
\begin{subequations}
    \begin{alignat}{2}
        && \Jcal(\v{\theta}^s) = & \Jcal(\tilde{\v{\theta}}^\star - \varepsilon^{\text{stepsize}} \tilde{\v{v}}) \label{eq:VF_expl1} \\
        && = & \sum_{i=1}^{N} \frac{1}{2} (\tilde{\v{\theta}}^\star - \varepsilon^{\text{stepsize}} \tilde{\v{v}})^\T 
        \hat{\v{H}}_{i,s} 
        (\tilde{\v{\theta}}^\star - \varepsilon^{\text{stepsize}} \tilde{\v{v}}) \notag \\
        &&&
        + \hat{\v{f}}_{i,s}^\T (\tilde{\v{\theta}}^\star - \varepsilon^{\text{stepsize}} \tilde{\v{v}}) 
        + \hat{c}_{i,s} \label{eq:VF_expl2} \\
        && = & \Jcal(\tilde{\v{\theta}}^\star) + \frac{1}{2} (\varepsilon^{\text{stepsize}} \tilde{\v{v}})^\T \hat{\v{H}}_{s} (\varepsilon^{\text{stepsize}} \tilde{\v{v}}) \notag \\
        &&& - \varepsilon^{\text{stepsize}} \tilde{\v{v}}^\T (\hat{\v{H}}_{s} \tilde{\v{\theta}}^\star + \hat{\v{f}}_{s}) \label{eq:VF_expl3} \\
        && \approx & \Jcal(\tilde{\v{\theta}}^\star) - \varepsilon^{\text{stepsize}} \tilde{\v{v}}^\T (\hat{\v{H}}_{s} \tilde{\v{\theta}}^\star + \hat{\v{f}}_{s}) \label{eq:VF_expl4}
    \end{alignat}
\end{subequations}
where $\hat{\v{H}}_{s} = \sum_{i=1}^N \hat{\v{H}}_{i,s}$, and $\hat{\v{f}}_{s} = \sum_{i=1}^N \hat{\v{f}}_{i,s}$. When $\varepsilon^{\text{stepsize}}$ is sufficiently small, the second term in \eqref{eq:VF_expl3} can be cancelled. Note that \eqref{eq:VF_expl4} comprises of two terms, $\Jcal(\tilde{\v{\theta}}^\star)$ indicates current optimal values. In the second term, $\hat{\v{H}}_{s} \tilde{\v{\theta}}^\star + \hat{\v{f}}_{s}$ indicates the gradient of current critical region, $\tilde{\v{v}}$ is the subgradient induced by explored critical regions in \eqref{eq:CRE.LSQ}.

Under the above notations, for every CSR conducted in RCDCRE, we have the following proposition
\begin{proposition}
    (see \cite[Appendix B]{guoRobustTieLineScheduling2018}):
    \label{prop.rcdcre.rot}
    In RCDCRE, suppose rotation martix $\v{R}$ grenerated from \eqref{eq:rcdcre.Givens} has been applied to rotate the coordinate system. Then at least a new critical region of $\tilde{\v{\theta}}^\star$ or a better parameter $\v{\theta}^s$ can be found for block coordinate descent in the new coordinate system.
\end{proposition}

\begin{proof}
    We first give a geometric intepratation of subgradient calculation in \eqref{eq:CRE.LSQ}, shown in Figure \ref{fig:subdifferential}. In \cite[Proposition 2.5]{larssonConditionalSubgradientOptimization1996}, it is referred to as conditional steepest descent direction.
    \begin{figure}[htbp]
        \centering   
        \vspace{-7pt}
        \includegraphics[width=0.24\textwidth]{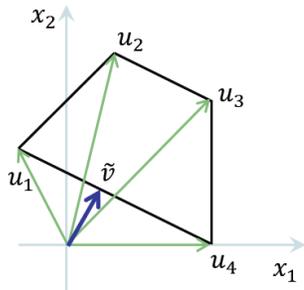}  
        \vspace{-7pt}
        \caption{An illustration to \eqref{eq:CRE.LSQ}.}
        \label{fig:subdifferential}
    \end{figure}   

    We let $\tilde{\v{v}}$ be solution to \eqref{eq:CRE.LSQ} with $\partial f(\v{x}^\star)$ and $\Ncal_{\v{\Xset}}^\star$. Here, $\partial f(\v{x}^\star)$ is the convex hull of subgradients $\v{u}_i, i = 1, \ldots, 4$, and $\Ncal_{\v{\Xset}}^\star$ is an empty set. 
    As $\tilde{\v{v}}$ is the negative of the shortest conditional subgradient to the conditional subdifferential $\partial f(\v{x}^\star) + \Ncal_{\v{\Xset}}^\star$.
    Clearly, if the conditional subdifferential does not contain the origin, then $\tilde{\v{v}} \neq \v{0}$ and $\tilde{\v{v}}^\T \v{u} > 0$ for any $\v{u} \in \partial f(\v{x}^\star) + \Ncal_{\v{\Xset}}^\star$.

    We distinguish two cases, if $(\hat{\v{H}}_{s} \tilde{\v{\theta}}^\star + \hat{\v{f}}_{s}) \in \partial \Jcal^\star$, we have
    \begin{alignat*}{2}
        \tilde{\v{v}}^\T (\hat{\v{H}}_{s} \tilde{\v{\theta}}^\star + \hat{\v{f}}_{s}) > 0
    \end{alignat*}
    Hence, by \eqref{eq:VF_expl4}, we have $\Jcal(\v{\theta}^s) < \Jcal(\tilde{\v{\theta}}^\star)$, which means $\tilde{\v{v}}$ is a descending direction.
    Now, suppose $\tilde{\v{v}}^\T (\hat{\v{H}}_{s} \tilde{\v{\theta}}^\star + \hat{\v{f}}_{s}) \leq 0$, implying $\Jcal(\v{\theta}^s) > \Jcal(\tilde{\v{\theta}}^\star)$. Although $\tilde{\v{v}}$ is not a descending direction, we can conclude that $(\hat{\v{H}}_{s} \tilde{\v{\theta}}^\star + \hat{\v{f}}_{s}) \notin \partial \Jcal^\star$ and we have explored a new critical region at $\tilde{\v{\theta}}^\star$. There are also cases that we can both find a better $\v{\theta}^s$ and a new critical region, we omit the similar analysis for brevity.
\end{proof}

\begin{corollary}
    (see \cite[Theorem 1]{guoRobustTieLineScheduling2018}):
    \label{corollary.rcdcre.converge}
    RCDCRE terminates after finitely many steps, when $\| \v{v} \| = 0$.
\end{corollary}

\begin{proof}
    As $\partial \Jcal^\star$ is a subdifferential at $\v{\theta}^\star$, it is sufficient to prove optimality when $\v{v} = \v{0}$, per first-order non-differentiable optimality condition $\v{0} \in \partial \Jcal^\star + \Ncal_{\v{\Theta}}^\star$. 
    We proof corollary \ref{corollary.rcdcre.converge} in two steps. Suppose $\v{\theta}^\star$ is non-optimal, then by proposition \ref{prop.rcdcre.rot}, a descending direction will be identified after finite CSRs as there are finite critical regions according to lemma \ref{lemma.mpP}. Hence, we can apply BCD as soon as a better $\v{\theta}^s$ has been found. When $\v{\theta}^\star$ becomes the optimal solutions, also by proposition \ref{prop.rcdcre.rot} and lemma \ref{lemma.mpP}, at least one new critical region will be identified for each CSR. Hence, RCDCRE converge in finite iterations.
\end{proof}

\subsubsection{Proof of \ref{proof:bigM})} 
We first show that the optimal solutions are identical for \eqref{eq:mpP} and \eqref{eq:mpP_bigM} when $\v{\theta}^s$ is feasible. Suppose the active set partition of \eqref{eq:mpP_bigM} is given by $\{\Acal, \bar{\Acal}\}$, $\{\Ical, \bar{\Ical}\}$, where $\Acal$, $\Ical$ are the partition of \eqref{eq:mpP}, by KKT conditions we have
\begin{subequations}
    \label{eq:mpP_bigM_KKT}
    \begin{alignat}{2}
        &&& \v{H} \v{x} + \v{f} + 
        \v{A}_\Acal^\T \v{\lambda}_\Acal = \v{0}, \label{eq:mpP_bigM_KKT1} \\
        &&& M \bone + \v{\lambda} - \v{\mu} = \v{0}, \label{eq:mpP_bigM_KKT2} \\        
        &&& \v{A}_\Acal \v{x} - \v{s}_{\Acal} = \v{b}_\Acal + \v{C}_\Acal \v{\theta},\label{eq:mpP_bigM_KK3} \\
        &&& \v{s}_{\bar{\Acal}} = \v{0}, \quad \v{\mu}_{\bar{\Ical}} = \v{0}, \quad \v{\lambda}_{\Ical} = \v{0}, \label{eq:mpP_bigM_KK4} \\
        &&& \v{A}_\Ical \v{x} - \v{s}_{\Ical} \leq \v{b}_\Ical + \v{C}_\Ical \v{\theta},\label{eq:mpP_bigM_KK5} \\
        &&& \v{s}_{\bar{\Ical}} \geq \v{0}, \quad \v{\mu}_{\bar{\Acal}} \geq \v{0}, \quad \v{\lambda}_{\Acal} \geq \v{0}. \label{eq:mpP_bigM_KK6}
    \end{alignat}
\end{subequations} 
Here, $\v{\lambda}$ is the multiplies to the initial constraints \eqref{eq:mpP.cons}, $\v{\mu}$ is the multiplies to the slack variables $\v{s}$.
As $M$ is sufficiently large, $\v{s}^\star = \v{0}$ when $\v{\theta}^s$ is feasible to \eqref{eq:mpP}, and $\v{\mu} = M \bone \geq \v{0}$. Equation \eqref{eq:mpP_bigM_KKT} reduces to 
\begin{alignat*}{2}
    &&& \v{H} \v{x} + \v{f} + 
    \v{A}_\Acal^\T \v{\lambda}_\Acal = \v{0}, \\    
    &&& \v{A}_\Acal \v{x} = \v{b}_\Acal + \v{C}_\Acal \v{\theta}, \quad \v{\lambda}_{\Acal} \geq \v{0}, \\
    &&& \v{A}_\Ical \v{x} \leq \v{b}_\Ical + \v{C}_\Ical \v{\theta}, \quad \v{\lambda}_{\Ical} = \v{0}. 
\end{alignat*}
which is exactly the KKT conditions of \eqref{eq:mpP}. Hence, the value function and critical regions of \eqref{eq:mpP} and \eqref{eq:mpP_bigM} are all the same for any feasible $\v{\theta}^s$. 
When $\v{\theta}^s$ is infeasible to \eqref{eq:mpP}, the problem can still be viewed as a feasible mp-LP/QP of \eqref{eq:mpP_bigM}. According to lemma \ref{lemma.mpP}, the value function is a continuous piecewise convex function, and there are finite critical region partitions. Clearly, the value function is continuous on the boundary of feasible and infeasible regions of \eqref{eq:mpP}. Hence, the effect of big$M$-penalty is merely to extend feasible space $\v{\Theta}^\star$.
By corollary \ref{corollary.rcdcre.converge}, we can conclude RCDCRE also converges to optimal solutions under big$M$-penalty formulation.

\subsubsection{Proof of \ref{proof:L1})}
The last concern of RCDCRE is the $\ell_1$-penalty formulation to the convergence of $\v{\theta}$ outside $\v{\Theta}$, which is summarized in the following lemma.
\begin{lemma}
    \label{lemma.L1Reg}
    (see \cite[5.16]{boydConvexOptimization2004})
    Suppose $\v{\nu}^\star$ is the optimal Lagrange dual solutions of \eqref{eq:rcdcre_mp_L1Reg}. If $\sigma > \bone^T \v{\nu}^\star$, the any solution to $\ell_1$-penalty problem \eqref{eq:rcdcre_mp_L1Reg} is also an optimal solution of the original problem.
\end{lemma}

So after several iterations when $\v{\Theta} \cap_{i=1}^N \text{CR}^\star_{i,s} \neq \emptyset$, the algorithm can return the same optimal solutions. The above arguments combines to prove Theorem \ref{thm.RCDCRE}.

\end{document}